\documentclass[11pt]{amsart}

\newtheorem{theorem}{Theorem}[section]

\theoremstyle{definition}

\theoremstyle{remark}

\numberwithin{equation}{section}

\newcommand\lam{\lambda}
\newcommand\na{\nabla}

\newcommand\Del{\Delta}

\newcommand\Rc{\textup{Rc}}

\newcommand\ppt{\frac{\partial}{\partial t}}
\newcommand\ddt{\frac{d}{dt}}

\begin{document}

\title[Monotonic Quantities]{A Class of Monotonic Quantities along the Ricci Flow}

\author[Jun Ling]{Jun \underline{LING}}
\address{Department of Mathematics, Utah Valley State College, Orem, Utah 84058}
\email{lingju@uvsc.edu}
\thanks{The author thank the Mathematical Sciences Research Institute at Berkeley for its
hospitality for program `The Geometric Evolution Equations and
Related Topics' and thank National Science Foundation for the
support offered.}


\subjclass[2000]{Primary 53C21, 53C44; Secondary 58J35, 35P99}

\date{October 12, 2007}


\keywords{Monotonicity, Ricci flow}

\begin{abstract}
We construct a class of monotonic quantities along the normalized
Ricci flow on closed $n$-dimensional manifolds.
\end{abstract}

\maketitle

\section{Introduction}\label{sec-intro}
The invariants and properties preserved along the Ricci flow often
play key roles in the study of geometry and topology of manifolds.
Among many of such examples are, positivity of scalar curvature
\cite{hamilton} in the study of three-manifolds with positive
Ricci curvature, positivity of curvature operators in the study of
high dimensional manifolds \cite{hamilton2}, Hamilton's entropy
and its monotonicity \cite{hamilton3} in the study of
two-dimensional manifolds, and Perelman's entropy and the
monotonicity of the the first eigenvalue of $-4\Del+R$ along the
Ricci flow \cite{perelman} in the study of high dimensional
manifolds, and etc. There has been increasing attention on the
last two examples. The monotonicity of the Hamilton's entropy was
a key to get an upper bound of scalar curvature $R$. With such a
bound and the Harnack inequality of LYH type, Hamilton
\cite{hamilton3} was able to get a crucial lower of $R$ and prove
the exponential convergence of the the normalized Ricci flow for
surfaces with positive Euler characteristic. With his entropy and
the monotonicity of the first eigenvalue, Perelman \cite{perelman}
was able to rule out nontrivial steady or expanding breathers on
compact manifolds. Hamilton and Perelman's work on the
monotonicity stimulated the research on the topic. Many studies on
the topics appeared. For the Laplacian operator, Li Ma \cite{mal}
studied monotonicity of the first eigenvalue on a domain $D$ in
the manifold $M$ with Dirichlet boundary condition, along the
unnormalized Ricci flow. The author \cite{lingj} recently studied
the first nonzero eigenvalue of Laplacian under the normalized
Ricci flow and gave a Faber-Krahn type of comparison theorem and a
sharp bound. In \cite{lingj2}, the author studied some asymptotic
behavior of the first nonzero eigenvalue of the Lalacian along the
normalized Ricci flow and gave a direct short proof for an
asymptotic upper limit estimate.

In this paper, we construct a class of monotonic quantities along
the normalized Ricci flow. We show that though the eigenvalues
themselves are not monotonic along the normalized Ricci flow in
general, the appropriate multiples are. We first introduce some
notations we are going to use in this paper.

We let $M$ be a closed $n$-dimensional manifold, $g(t)$ the
solution to the normalized Ricci flow equation
\begin{equation}\label{nrf}
\ppt g=-2\Rc+\frac{2r}{n}g\qquad \textup{for}\quad 0\leq
t<T\leq\infty.
\end{equation}
Let $\Rc$ be the Ricci tensor, $R$ the scalar curvature, $d\mu$
the volume element, $\Del$ the Laplacian, of the Riemannian
manifold $(M,g(t))$, respectively. Let
\[
r=r(t)=\int_MRd\mu\Big/\int_Md\mu,\quad
\sigma(t)=\int_0^tr(\tau)d\tau,
\]
\[
\rho_0=\min_MR|_{t=0}\qquad \delta_0=\max_MR|_{t=0}.
\]
We drop the integral domain $M$ in integrals sometimes.

We present the non-decreasing quantities in the next section and
increasing ones in the last section.

\section{Non-decreasing quantities}\label{sec-pe}

We have the following results on non-decreasing quantities along
the normalized Ricci flow.

\begin{theorem}\label{thm1+}
Let $g(t)$ be the solution to the normalized Ricci flow equation
(\ref{nrf}), $\lam=\lam(t)$ be any eigenvalue of the Laplacian of
the metric $g(t)$ on a closed $n$-dimensional manifold $M$. If the
Einstein tensor $E=:\Rc-\frac12Rg$ is non-negative, then the
quantity
\[
e^{\int_0^t[\frac2n r(\tau)-\varphi(\tau)]d\tau}\lam(t)
\]
is non-decreasing along the flow, where
\begin{equation}\label{varphi-def}
\varphi(t)=1\Big/\left\{e^{\frac2n\sigma(t)}\left(\frac{1}{\rho_0}
-\frac2n\int_0^te^{-\frac2n\sigma(\tau)}d\tau\right)\right\}.
\end{equation}
\end{theorem}

\begin{proof}
Let $u$ be an eigenfunction of the eigenvalue of the Laplacian,
\[
-\Del u =\lam u.
\]
Take derivatives with respect to $t$,
\[
-(\ppt \Del) u -\Del \ppt u =(\ddt\lam) u+ \lam\ppt u.
\]
Multiply the equation by $u$ and integrate,
\[
-\int u(\ppt \Del) u  -\int u\Del \ppt u = (\ddt \lam) \int u^2+
\lam\int u\ppt u.
\]
Noticing that
\[
-\int u\Del \ppt u=-\int\Del u\, \ppt u=\lam\int u\ppt u,
\]
we have
\begin{equation}\label{ddt-lam-1}
\begin{split}
&(\frac{d}{dt}\lam)\int u^2d\mu\\
&=-\int u(\ppt \Del) ud\mu
\\
&=-\int (2R^{ij}u\na_i\na_ju-\frac{2r}{n}u\Del u)d\mu\\
&=-\int 2R^{ij}u\na_i\na_jud\mu-\int\frac{2r}{n}\lam u^2d\mu,
\end{split}
\end{equation}
where in the second equality we used the equation
\begin{equation}\label{ddt-laplacian}
\ppt(\Del)=2R^{ij}\na_i\na_j-\frac{2r}{n}\Del.
\end{equation}
This equation is true is due to the following. By (\ref{nrf}), for
a smooth function $v$ in a local chart $\{x^i\}$ on $M$,
\[
\begin{split}
\ppt(\Del_{g(t)}v)&=\ppt(g^{ij}\na_i\na_jv)
\\
&=(\ppt
g^{ij})\na_i\na_jv+g^{ij}\big[\frac{\partial^2}{\partial_x^i\partial
x^j}\ppt v -\Gamma_{ij}^k\frac{\partial}{\partial x^k}\ppt
v-\ppt(\Gamma_{ij}^k)\frac{\partial}{\partial x^k}v\big]
\\
&=-g^{ik}g^{jl}(-2R_{kl}+\frac{2r}{n}g_{kl})\na_i\na_j
v+\Delta\ppt v-g^{ij}\ppt(\Gamma_{ij}^k)\frac{\partial}{\partial
x^k}v
\\
&=2R^{ij}\na_i\na_jv-\frac{2r}{n}\Del
v-g^{ij}\ppt(\Gamma_{ij}^k)\frac{\partial}{\partial
x^k}v+\Delta\ppt v,
\end{split}
\]
and at point $x$ and in the local normal chart about $x$,
\[
g^{ij}\ppt(\Gamma_{ij}^k)=\frac12g^{ij}g^{kl}[\frac{\partial}{\partial
x^i}g_{lj}+\frac{\partial}{\partial
x^j}g_{il}-\frac{\partial}{\partial x^l}g_{ij}]
\]
\[
=\frac12g^{kl}[g^{ij}\nabla_i\ppt g_{lj}+g^{ij}\nabla_j\ppt
g_{il}-g^{ij}\nabla_l\ppt g_{ij}]=0,
\]
where in the last equality we used (\ref{nrf}) and the contracted
second Bianchi identity. Therefore (\ref{ddt-laplacian}) holds.

Now the first term in (\ref{ddt-lam-1}) is
\[
-\int2R^{ij}u\na_i\na_jud\mu=\int(2\na_iR^{ij})u\na_jud\mu+\int2R^{ij}\na_iu\na_jud\mu.
\]
The contracted second Bianchi identity implies that
\[
\begin{split}
\int(2\na_iR^{ij})u\na_jud\mu &=\int u(\na^jR)\na_jud\mu \\
&=-\int Ru\Del ud\mu-\int R|\na u|^2d\mu \\
&=\lam\int R u^2d\mu -\int R|\na u|^2d\mu.
\end{split}
\]
The above two equations give
\[
-\int2R^{ij}u\na_i\na_jud\mu=\lam\int R u^2d\mu-\int R|\na
u|^2d\mu+\int2\Rc(\na u,\na u)d\mu.
\]
Therefore (\ref{ddt-lam-1}) becomes
\[
\begin{split}
&(\frac{d}{dt}\lam)\int u^2d\mu\\
 &=\lam\int R u^2 -\int R|\na
u|^2+\int2\Rc(\na u,\na u) -\frac{2r}{n}\lam\int u^2d\mu
\\
&=\lam\int \big[R-\frac{2}{n}r\big]u^2d\mu -\int R|\na u|^2d\mu
+\int2\Rc(\na u,\na u)d\mu,
\end{split}
\]
that is,
\begin{equation}\label{lam-evo-eq}
\begin{split}
\frac{d}{dt}\lam &=\frac{\int_M
\big[R-\frac{2}{n}r\big]u^2d\mu}{\int_M
u^2d\mu}\lam+\frac{2\int_M[\Rc(\na u,\na u)-\frac12R|\na
u|^2]d\mu}{\int_M u^2d\mu}.
\end{split}
\end{equation}

On the other hand, the evolution equation of $R$
 \[
\ppt R=\Del R + 2|\Rc|^2-\frac{2r}{n}R.
 \]
 and the inequality
\[
|\Rc|^2\geq \frac1nR^2
\]
imply that
\begin{equation}\label{R-evo-ineq}
\ppt R\geq\Del R + \frac2nR(R-r).
\end{equation}
Therefore by the maximum principle, we have
\[
R\geq\varphi(t),
\]
where $\varphi$, as defined in (\ref{varphi-def}), is the solution
to the ODE initial value problem
 \[
\ddt \varphi=\frac2n\varphi(\varphi-r), \quad \varphi(0)=\rho_0.
 \]
Taking this time-dependent lower bound of $R$ into
(\ref{lam-evo-eq}) and using the non-negativity of the Einstein
tensor, we get
\[
\ddt\lam\geq [\varphi(t)-\frac2n r(t)]\lam
\]
and
\[
\ddt
\left\{e^{\int_0^t[\frac2nr(\tau)-\varphi(\tau)]d\tau}\lam(t)\right\}\geq
0.
\]
\end{proof}

\begin{theorem}\label{thm2+}
Let $g(t)$ be a solution to the normalized Ricci flow (\ref{nrf})
on a closed two-dimensional manifold, $\lam=\lam(t)$ any
eigenvalue of the Laplacian $\Del$ of $(M,g(t))$. Then  the
quantity
\[
\left|\frac{\rho_0}{r_0}-\frac{\rho_0}{r_0}e^{r_0\,t}+e^{r_0\,t}\right|\lam(t)
\]
is non-decreasing if the Euler characteristic $\chi\not=0$; and
the quantity
\[
(1-\rho_0 t)\lam(t)
\]
is non-decreasing if the Euler characteristic $\chi=0$, along the
normalized Ricci flow $g(t)$.
\end{theorem}

\begin{proof}
By (\ref{nrf}) we have
\[
\frac{d}{dt}(d\mu) =(r-R)d\mu
\]
and
\[
\frac{d}{dt}\int_M  d\mu=\int_M(r-R)d\mu=0.
\]
So the volume of $M$ remains constant in $t$ along the flow.

Now in dimension $n=2$, the Gauss-Bonnet Theorem implies that
\begin{equation}\label{r-chi}
r=4\pi\chi\Big/\int_Md\mu.
\end{equation}
Therefore $r$ is a constant, $r\equiv r_0$.

Note that in dimension two the evolution of $R$ is
\[
\ppt R=\Del R + R(R-r).
\]
Therefore we have by the maximum principle
\[
R\geq \varphi(t),
\]
where $\varphi$ is the solution to the ODE initial value problem
\[
\ddt\varphi=\varphi(\varphi-r),\qquad
\varphi(0)=\rho_0=:\min_MR|_{t=0}.
\]
Note that
\[
\varphi(t)=\frac{r_0}{1-\left(1-\frac{r_0}{\rho_0}\right)e^{r_0t}}
\]
in the case of $\chi\not=0$;
\[
\varphi(t)=\frac{\rho_0}{1-\rho_0t}
\]
in the case $\chi=0$.

Taking the above time-dependent lower bound $\varphi(t)$ into
(\ref{lam-evo-eq}), noticing that the second term is zero since
$\Rc=\frac12Rg$ in dimension two, and integrating, we get the
theorem.
\end{proof}

\section{Non-increasing Quantities}\label{sec-ne}

We have the following results on non-increasing quantities along
the normalized Ricci flow.

\begin{theorem}\label{thm1-}
Let $g(t)$ be the solution to the normalized Ricci flow equation
(\ref{nrf}), $\lam=\lam(t)$ be any eigenvalue of the Laplacian of
the metric $g(t)$ on a closed $n$-dimensional manifold $M$. If the
Ricci tensor $\Rc$ is non-negative and the Einstein tensor
$E=\Rc-\frac12Rg$ is non-positive, that is,
\[
0\leq \Rc\leq \frac12Rg,
\]
then the quantity
\[
e^{\frac2n\sigma(t)-\int_0^t\psi(\tau)d\tau}\lam(t)
\]
is non-increasing along the normalized Ricci flow, where
\begin{equation}\label{psi-def}
\psi(t)=1\Big/\left\{e^{\frac2n\sigma(t)}
\left(\frac1{\delta_0}-2\int_0^te^{-\frac2n\sigma(\tau)}d\tau\right)\right\}.
\end{equation}
\end{theorem}

\begin{proof}
That $\Rc\geq 0$ implies $|\Rc|^2\leq R^2$. Taking this into the
evolution equation of $R$
\[
\ppt R=\Del R+2|\Rc|^2-\frac{2}{n}rR,
\]
we get
\[
\ppt R\leq \Del R+2R(R-\frac1nr).
\]
Compare the above inequality with the ODE
\[
\ddt \psi=2\psi(\psi-\frac1nr).
\]
Let $\psi(t)$ be the solution of the ODE with the initial value
\[
\psi(0)=\delta_0=:\max_MR|_{t=0}.
\]
It is easy to see that $\psi$ is the function defined in
(\ref{psi-def}). The maximum principle implies that
\[
R\leq \psi(t).
\]
Taking this into (\ref{lam-evo-eq}), and noticing the
non-positivity of the Einstein tensor, we get
\[
\frac{d}{dt}\lam \leq\frac{\int
\big[R-\frac{2}{n}r\big]u^2d\mu}{\int u^2d\mu}\,\lam\leq
(\psi(t)-\frac2nr)\lam.
\]
and
\[
\ddt (e^{\frac2n\sigma(t)-\int_0^t\psi(\tau)d\tau}\lam(t))\leq 0.
\]
\end{proof}

\begin{theorem}\label{thm2-}
Let $g(t)$ be a solution to the normalized Ricci flow (\ref{nrf})
on a closed two-dimensional manifold, $\lam=\lam(t)$ any
eigenvalue of the Laplacian $\Del$ of $(M,g(t))$. Then the
quantity
\[
\left|\frac{\delta_0}{r_0}-\frac{\delta_0}{r_0}e^{r_0\,t}+e^{r_0\,t}\right|\lam(t)
\]
is non-increasing if the Euler characteristic $\chi\not=0$; and
the quantity
\[
(1-\delta_0 t)\lam(t)
\]
is non-decreasing if the Euler characteristic $\chi=0$, along the
normalized Ricci flow $g(t)$.
\end{theorem}

\begin{proof}
The same argument as in the proof of Theorem \ref{thm2+} shows
that $r$ is a constant independent of $t$, $r\equiv r_0$.

By the evolution equation of $R$ in dimension two
\[
\ppt R=\Del R+ R(R-r)
\]
and the maximum principle, we have
\[
R\leq \psi(t),
\]
where  $\psi(t)$ is the solution to the ODE initial value problem
\[
\ddt \psi=\psi(\psi-r),\qquad \psi(0)=\delta_0=:\max_MR|_{t=0}.
\]
It is easy to know that
\[
\psi(t)=\frac{r_0}{1-\left(1-\frac{r_0}{\delta_0}\right)e^{r_0t}}
\]
 in the case of $\chi\not=0$, and
\[
\psi(t)=\frac{\delta_0}{1-\delta_0t},
\]
in the case $\chi=0$.

Taking this into (\ref{lam-evo-eq}) and noticing $\Rc=\frac12Rg$
in dimension two, we get the theorem.

\end{proof}

\begin{center}
\textbf{Acknowledgement}
\end{center}
The author is grateful to Professor Yieh-Hei Wan for his constant
support, valuable suggestions and helps,  to Professor Bennett
Chow for his many enlightening talks, papers and books,
suggestions and helps.

\end{document}